\documentclass{amsart}

\usepackage{amsmath}
\usepackage{amssymb}
\usepackage[all]{xy}
\usepackage{picinpar}
\usepackage{palatino}

\def\eu{\mathfrak}
\def\ma{\mathbb}

\def\finite#1{{\ma F}_{q^{#1}}}

\def\lam#1{k(\Lambda_#1)}

\def\G#1{\big(R_T/(#1)\big)^{\ast}}
\def\g#1{{#1_{{\eu g}{\eu e}}}}
\def\p{{\mathcal P}_{\infty}}
\def\f{{\ma F}_q^{\ast}}
\def\F{{\ma F}_q}
\def\producto{{\eu p}_1^{e_1}\cdots {\eu p}_r^{e_r}}
\def\pinfty{{\eu p}_{1,{\infty}}^{e_{1,{\infty}}}
\cdots {\eu p}_{r_{\infty},\infty}^{e_{r_{\infty},\infty}}}
\def\P{{\mathcal P}}

\def\fin{\hfill\qed\bigskip}

\newcommand{\Irr}{\operatorname{Irr}}
\newcommand{\Gal}{\operatorname{Gal}}
\newcommand{\con}{\operatorname{con}}

\newcommand{\lcm}{\operatorname{lcm}}

\newcounter{bean}

\def\l{
\begin{list}
{\rm{(\alph{bean}).-}}{\usecounter{bean}
\setlength{\labelwidth}{0.8in}
\setlength{\labelsep}{0.3cm}
\setlength{\leftmargin}{1cm}}}

\numberwithin{equation}{section}
\newtheorem{theorem}{Theorem}[section]
\newtheorem{proposition}[theorem]{Proposition}
\newtheorem{lemma}[theorem]{Lemma}
\newtheorem{example}[theorem]{Example}
\newtheorem{remark}[theorem]{Remark}
\newtheorem{definition}[theorem]{Definition}

\title[Genus Fields of Congruence Function Fields]
{Genus Fields of Congruence Function Fields}

\author[M. Maldonado]{Myriam Maldonado--Ram\'irez}
\address{Departamento de Matem\'aticas\\
Escuela Superior de F\'isica y Matem\'aticas del I.P.N.}
\email{rosalia@esfm.ipn.mx}

\author[M. Rzedowski]{Martha Rzedowski--Calder\'on}
\address{Departamento de Control Autom\'atico\\
Centro de Investigaci\'on y de Estudios Avanzados del I.P.N.}
\email{mrzedowski@ctrl.cinvestav.mx}

\author[G. Villa]
{Gabriel Villa--Salvador}
\address{
Departamento de Control Autom\'atico\\
Centro de Investigaci\'on y de Estudios Avanzados del I.P.N.}
\email{gvillasalvador@gmail.com, gvilla@ctrl.cinvestav.mx}



\date{December 27th., 2015}

\begin{document}

\begin{abstract}

Let $k$ be a rational congruence function field and
consider an arbitrary finite separable
extension $K/k$. If for each
prime in $k$ ramified in $K$ we have that at least one
ramification index is not divided by the characteristic of $K$,
we find the genus field $\g K$, except for constants, 
of the extension $K/k$. In general, we describe the
genus field of a global function field.

\end{abstract}

\maketitle

\section{Introduction}\label{S1}

C. F. Gauss \cite{Gau1801} was the first in considering what
now is called the {\em genus field} of a quadratic
number field. H. Hasse \cite{Has51} introduced genus theory
for quadratic number fields describing the theory invented by
Gauss by means of class field theory. H. W. Leopoldt \cite{Leo53}
determined the genus field $\g K$ of an absolute abelian
number field $K$ generalizing the work of Hasse. Leopoldt
developed the theory using Dirichlet characters and relating them
with the arithmetic of $K$. A. Fr\"ohlich \cite{Fro59-1, Fro59-2} introduced
the concept of genus fields for nonabelian number fields. 
Fr\"ohlich defined the genus field $\g K$ of an arbitrary
finite number field $K/{\ma Q}$ as $\g K:= Kk^{\ast}$ where
$k^{\ast}$ is the maximal abelian number field such that $Kk^{\ast}/
K$ is unramified. We have that $k^{\ast}$ is the maximal abelian
number field contained in $\g K$. The degree $[\g K:K]$ is called
the {\em genus number} of $K$ and the Galois group $\Gal(
\g K/K)$ is called the {\em genus group} of $K$.

If $K_H$ denotes the Hilbert class field (HCF) of $K$, we have
$K\subseteq \g K\subseteq K_H$ and $\Gal(K_H/K)$ is
isomorphic to the class group $Cl_K$ of $K$. Then $\g K$
corresponds to a subgroup $G_K$ of $Cl_K$, that is,
$\Gal(\g K/K)\cong Cl_K/G_K$. The subgroup $G_K$ is called
the {\em principal genus} of $K$ and $|Cl_K/G_K|$ is equal to
the genus number of $K$.

M. Ishida \cite{Ish76} described the genus field $\g K$ of any
finite extension of ${\ma Q}$, allowing ramification at the
infinite primes. Given a number field $K$, Ishida found an
abelian number field $k_1^{\ast}$ and described another number
field $k_2^{\ast}$ such that $k^{\ast}=k_1^{\ast}k_2^{\ast}$ and
$k_1^{\ast}\cap k_2^{\ast}={\ma Q}$. The field $k_1^{\ast}$ was
constructed by means of the finite primes $p$ such that at least
one ramification index of the decomposition of $p$ in $K$ is not
divisible by $p$. In other words, by those primes $p$ such that
at least one prime in $K$ above $p$ is tamely ramified.

We are interested in genus fields in the context of congruence
(global) function fields. In this case there is no proper notion of
Hilbert class field since all the constant field extensions are abelian
and unramified. In fact, if the class number of a congruence function
field $K$ is $h_K$ then there are exactly $h:=h_K$ abelian
extensions $K_1,\ldots, K_h$ of $K$ such that $K_i/K$ are maximal
unramified with exact field of constants of each $K_i$ the same
as the one of $K$, ${\ma F}_q$, the finite field of $q$ elements
and $\Gal(K_i/K)\cong Cl_{K,0}$ the group of classes of
divisors of degree zero
(\cite[Chapter 8, page 79]{ArTa67}).

M. Rosen \cite{Ros87} gave a definition of Hilbert class fields of $K$,
fixing a nonempty finite set $S_{\infty}$ of prime divisors of $K$.
Using Rosen's definition of HCF is possible to give a proper concept
of genus fields along the lines of number fields. In the literature
there have been different definitions of genus fields according to 
different HCF definitions. R. Clement \cite{Cle92} found a narrow genus
field of a cyclic extension of $k=\F(T)$ of prime degree dividing
$q-1$. She used the concept of HCF similar to that of a quadratic
number field $K$: it is the finite abelian extension of $K$ such
that the prime ideals of the ring of integers ${\mathcal O}_K$
of $K$ splitting are precisely the principal ideals generated by 
an element of positive norm. S. Bae and J. K. Koo \cite{BaKo96}
were able to generalize the results of Clement with the methods
developed by Fr\"ohlich \cite{Fro83}. They defined the genus
field for general global function fields and developed the analogue
of the classical genus theory. 

B. Angl\`es and J.-F. Jaulent
\cite{AnJa2000} used narrow $S$--class groups to establish
the fundamental results for genus theory of finite extensions of
global fields, where $S$ is a finite nonempty set of places.
G. Peng \cite{Pen2003} explicitly described the genus theory for
Kummer function fields for a prime number $l$
based on the global function field analogue of P. E.
Conner and J. Hurrelbrink exact hexagon.

C. Wittman \cite{Wit2007} extended Peng's results to the
case $l\nmid q(q-1)$ and used his results to study the $l$--part
of the ideal class groups of cyclic extensions of degree $l$.
S. Hu and Y. Li \cite{HuLi2010} described explicitly the genus
field of an Artin--Schreier extension of $k=\F(T)$.

In \cite{MaRzVi2013, MaRzVi2015} we
developed a theory of genus fields of
congruence function fields using Rosen's definition of HCF.
The methods we used there 
are based on the ideas of Leopoldt using Dirichlet
characters and a general description of $\g K$ in terms of
Dirichlet characters. The genus field 
$\g K$ was obtained for an abelian extension $K$ of
$k$. The method can be used to give $\g K$ explicitly when
$K/k$ is a cyclic extension of prime degree $l\mid q-1$ (Kummer)
or $l=p$ where $p$ is the characteristic (Artin--Schreier) and
also when $K/k$ is a $p$--cyclic extension (Witt). Later
on, the method was used in \cite{BaRzVi2013} to describe
$\g K$ explicitly when $K/k$ is a cyclic extension of degree
$l^n$, where $l$ is a prime number and $l^n\mid q-1$.

In this paper we consider a congruence function field $K$ with
exact field of constants a finite extension of $\F$ and such that $K$
is a finite separable extension of $k=
{\ma F}_q(T)$. We use 
Rosen's definition of HCF
of $K$ to find $\g K$ as $\g K=Kk^{\ast}$, where $k^{\ast}$ is the
composite of two fields $k_1^{\ast}$ and $k_2^{\ast}$.
Following Ishida's methods, we prove that $k_1^{\ast}$ is 
contained in the composite of constants with
fields $F_{\P}$ where
$k\subseteq  F_{\P} \subseteq \lam P$, $\P$
is fully ramified in $F_{\P}/k$, $\P$ running in the set of finite
primes ramified in $K/k$. Here $P$ is the monic irreducible
polynomial in $T$ associated to $\P$ and $\lam P$
is the $P$--th cyclotomic function field. The field $k_2^{\ast}$
encodes the wild ramification of the extension $k^{\ast}/k$.
The main difficulty handling $k_1^{\ast}$
is the decomposition of the infinite primes.

\section{Notation and basic results on cyclotomic function fields}\label{S2}

The results on function fields and cyclotomic function fields 
we need in this
paper may be consulted in \cite{Vil2006}. Let $\F$ be the finite
field of $q$ elements and of characteristic $p$. Let $k=\F(T)$ be a fixed
rational function field. Let $K$ be a congruence (global) function field with
exact field of constants a finite extension of $\F$ and
such that $K/k$ is a finite separable
extension. Let $R_T=\F[T]$ be the ring
of polynomials and $R_T^+$ denotes the subset of monic irreducible
polynomials. The pole of $T$ in $k$ will be denoted by $\p$. We
say that $\p$ is the {\em infinite prime} of $k$. 

For $M\in R_T\setminus
\{0\}$, $\Lambda_M$ denotes the $M$ torsion of the Carlitz module 
and $\lam M$ denotes the $M$--th cyclotomic function field.
We have that $\lam M/k$ is an abelian extension and $\Gal(\lam M/k)
\cong \G M$. For any $M\in R_T\setminus \{0\}$ the prime $\p$ has
ramification index $q-1$ and
decomposes in $\Phi(M)/(q-1)$ primes of degree one in $\lam M$.
The inertia group of $\p$ in $\lam M/k$ is identified with $\f
\subseteq \G M$. The fixed field $\lam M^{\F^{\ast}}=\lam M^+$
is the {\em maximal real subfield} of $\lam M$. 

By a geometric extension we mean an extension without new
constants.

The definition we use for the Hilbert class field is the one given by
Rosen. That is:

\begin{definition}\label{D2.1}{\rm{
Given a congruence function field $K$, the {\em Hilbert class field}
$K_H$ of $K$ is defined as the maximal unramified abelian extension
of $K$ such that all the primes in $K$ above $\p$ decompose fully.
}}
\end{definition}

With this definition of Hilbert class field, the definition of the genus
field of $K/k$ is given as follows.

\begin{definition}\label{D2.2}{\rm{
Let $K$ be a finite separable 
extension of $k$. The {\em genus field} $\g K$ of 
$K$ is the maximal extension of $K$ contained in $K_H$ that is the
composite of $K$ and an abelian extension of $k$. Equivalently, 
$\g K=K k^{\ast}$, where $k^{\ast}$ is the maximal abelian extension of
$k$ contained in $K_H$.
}}
\end{definition}

Note that it is possible to have $\g K=K E$ for several different
subfields $E\subsetneqq k^{\ast}$. We are interested in $k^{\ast}$
itself. In particular we have $k^{\ast}=\g {k^{\ast}}$ and since $k^{\ast}/k$
is abelian, the description of such $k^{\ast}$ may be found in
\cite{MaRzVi2013}.

The set of prime divisors in $k$ will be denoted by ${\ma P}_k$ and
let ${\ma P}_k^{\ast}:={\ma P}_k\setminus \{\p\}$ be the set of
finite primes of $k$.

The conorm map
from a field $E$ to a field $F$ will be denoted by $\con_{E/F}$. In an extension
$F/E$, $e(F|E)$ denotes the ramification index of a prime in $F$ above
one in $E$. If the primes are ${\eu P}$ and ${\eu p}$ we also write
$e(F|E)=e({\eu P}|{\eu p})=e_{F/E}({\eu P}|{\eu p})$.
The symbol $d_E({\eu p})$
denotes the degree of ${\eu p}$ for a prime ${\eu p}$ in $E$.
If $\P$ is a prime in $k$, its degree will be denoted by $d_{\P}$.

For any finite extension $E/k$ and any $m\in{\ma N}$ we denote
the extension of constants $E{\ma F}_{q^m}$ by $E_m$.

Let $\P_1,\ldots, \P_s,\P_{s+1},\ldots, \P_t$ 
be the finite primes in $k$ ramified in $K$. Let
$P_i\in R_T^+$ be such that the divisor $(P_i)_k$ is 
$(P_i)_k=\frac{\P_i}{\p^{\deg P_i}}$ for $1\leq i\leq t$.
For a prime $\P\in{\ma P}_k$, if $\con_{k/K}\P=\producto$,
we denote 
\begin{gather}\label{Eq0.1}
e_{\P}=\gcd(e_1,\ldots,e_r)= p^{u_{\P}} e_{\P}^{(0)}, \quad u_{\P}\geq
0,\quad \gcd(p,e^{(0)}_{\P})=1.
\end{gather}
We assume that
$p\nmid e_{\P_i}$ for $1\leq i\leq s$ and $p\mid e_{\P_j}$ for
$s+1\leq j\leq t$. That is, $u_{\P_i}=0$ for
$1\leq i\leq s$ and $u_{\P_j}
\geq 1$ for $s+1\leq j\leq t$.

One of the main tools used in this paper is the following result.

\begin{theorem}[Abhyankar's Lemma]\label{T2.3}
Let $F/E$ be a finite separable extension of function fields. Suppose that
$F=E_1E_2$ with $E\subseteq E_i\subseteq F$. Let ${\P}$ be a prime
of $E$ and ${\eu P}$ a prime in $F$ above ${\P}$. Let 
${\eu p}_i={\eu P}
\cap E_i$ for $i=1,2$. If at least one of the extensions $E_i/E$ is tamely
ramified at ${\eu p}_i$, then
\[
e_{F/E}({\eu P}|{\P})=\lcm[e_{E_1/E}({\eu p}_1|{\P}), e_{E_2/E}
({\eu p}_2|{\P})].
\]
\end{theorem}

\proof \cite[Theorem 12.4.4]{Vil2006}. $\fin$

We also recall two results.

\begin{proposition}\label{PalestineP4.1} Let $L/k$ be a finite abelian
extension, $P\in R_T^+$   and  $d:=\deg P$. Assume $P$
is tamely ramified in $L/k$. If $e$ denotes
the ramification index of $P$ in $L/k$, we have
$e\mid q^d-1$.
\end{proposition}

\proof \cite[Proposition 4.1]{SaRzVi2013}. $\fin$

\begin{proposition}\label{PalestineP4.2}
Let $L/k$ be an abelian extension where at most
a prime divisor ${\eu p}_0$ of degree $1$ is
ramified and the extension is tamely ramified. Then
$L/k$ is a constant extension.
\end{proposition}

\proof \cite[Proposition 4.2]{SaRzVi2013}. $\fin$

With respect to the genus field of a finite abelian
extension of $k$, we have the following results (see \cite[Theorem
4.2]{MaRzVi2013, MaRzVi2015}).

\begin{theorem}\label{T2.6}
Let $K/{\ma F}_q$ be a finite abelian
extension of $k$ where $\p$ is tamely ramified. Let $N\in R_T$
and $m\in{\ma N}$ be such that $K\subseteq {\lam N}{\ma F}_{q^m}$.
Let $E_{{\eu {ge}}}$ be the genus field of $E:=k(\Lambda_N)\cap
K{\ma F}_{q^m}$. 
Let $H_1$ be the subgroup that corresponds
to the decomposition group of the infinite primes of $K$ in
$\g E K/K$ under the Galois correspondence. Then
the genus field of $K$ is 
$$
K_{{\eu {ge}}}= \g {E^{H_1}}K.
\eqno{\Box}
$$
\end{theorem}

\begin{remark}\label{R2.7} {\rm{In Theorem \ref {T2.6}, it was assumed
originally that
$K/k$ is a geometric extension. This is not necessary; the same proof
that works for the geometric case works for any
finite abelian extension $K/k$.
}}
\end{remark}

\section{General case}\label{S4}

We consider a finite separable extension $K/k$ and let
$\g K=Kk^{\ast}$.
First we prove some general results.

\begin{proposition}\label{P3.1} 
Let $E/k$ be a finite tamely ramified  abelian extension.
For a finite prime ${\mathcal P}\in {\ma P}_k$, let 
$\con_{k/K}{\mathcal P}=\producto$ and
let $e^{\ast}_{\P}$ be the ramification index of $\P$ in $E/k$. 
Then $KE/K$ is
unramified at the primes above $\P$ if and only if 
$e^{\ast}_{\P}\mid e_{\P}$, where $e_{\P}$ is given
by {\rm{(\ref{Eq0.1})}}.
\end{proposition}

\proof Let ${\eu P}$ be a prime in $KE$ above
 $\P={\eu P}\cap k$. Thus
${\eu P}\cap K={\eu p}_i$ for some $i$.
Then, from Abhyankar's Lemma, we
have 
\[
e({\eu P}|\P)=\lcm[e({\eu p}_i|\P),e({\eu P}\cap E|\P)]=
\lcm [e_i,e^{\ast}_{\P}]=e({\eu P}|{\eu p}_i)e({\eu p}_i|\P)=
e({\eu P}|{\eu p}_i) e_i.
\]

Therefore ${\eu P}$ is unramified in $KE/K
\iff e({\eu P}|{\eu p}_i)=1 \iff \lcm[e_i,e^{\ast}_{\P}]=
e_i\iff e_{\P}^{\ast}\mid e_i$. The result follows. $\fin$

Consider the conorm of the infinite prime of $k$:
\begin{gather}\label{Eq3.2'}
\con_{k/K} \p=\pinfty. 
\end{gather}
Let $t_i$ be the degree of ${\eu p}_{i,\infty}$ and
\begin{gather}\label{Eq3.2''}
t_0:=\gcd (t_1,\ldots,t_{r_{\infty}}).
\end{gather} 

\begin{proposition}\label{P3.2} The field of constants of $\g K$ is
${\ma F}_{q^{t_0}}$.
\end{proposition}

\proof Consider the extension of constants $K{\ma F}_{q^m}/K$ with
$m\geq 1$. We have that ${\eu p}_{i,\infty}$ splits into $\gcd(t_i,m)$ factors
(\cite[Theorem 6.2.1]{Vil2006}). Therefore ${\eu p}_{i,\infty}$ 
decomposes fully
in $K{\ma F}_m \iff \gcd(t_i,m)=m\iff m\mid t_i$.
Thus the infinite primes of $K$ decompose 
fully in $K{\ma F}_m \iff m\mid t_0$. It
follows that ${\ma F}_{q^{t_0}}$ is the field of constants of $\g K$. $\fin$

\begin{proposition}\label{P3.3} Let $\P$ 
be a prime divisor of $k$ of degree 
$d$ such that $\P\neq \p$ and
$p\nmid e_{\P}$. Let $\g K=Kk^{\ast}$ and 
let $e_{\P}^{\ast}$ be the 
ramification index of $\P$ in $k^{\ast}/k$. 
Then $\gcd\big(e_{\P},
\frac{q^d-1}{q-1}\big)\mid e_{\P}^{\ast}$ 
and $e_{\P}^{\ast}\mid e_{\P}$.
\end{proposition}

\proof From Proposition \ref{P3.1} we have $e^{\ast}_{\P} \mid 
\gcd (e_1,\ldots,e_r)=e_{\P}$. Furthermore $\P$ is fully ramified
in $\lam P/k$, where $(P)_k=\frac{\P}{\p^{\deg P}}$ and
$\p$ decomposes fully in $\lam P^{\f}/k$. The degree
of the extension $\lam P^{\f}/k$
is $(q^d-1)/(q-1)$. Let $S$ be the subfield $k\subseteq
S\subseteq \lam P$ of degree $\gcd\big(e_{\P},\frac{q^d-1}{q-1}\big)$.
Then by Proposition \ref{P3.1} we have that
 $S$ satisfies that $KS/K$ is unramified and the infinite primes in $K$
decompose fully in $KS/K$ since $\p$ decomposes fully in $S/k$. 
Therefore $KS\subseteq Kk^{\ast}$, $S\subseteq k^{\ast}$ and
$\gcd\big(e_{\P},\frac{q^d-1}{q-1}\big)\mid e_{\P}^{\ast}$. $\fin$

Let $G:=\Gal(k^{\ast}/k)$ and let $G_p$ be the $p$--Sylow subgroup
of $G$. Then $G=G_0\times G_p$ with $p\nmid |G_0|$.
Therefore we have the decomposition
\begin{gather}\label{Eq3.1''}
k^{\ast}=k_1^{\ast}k_2^{\ast},\quad k_1^{\ast}\cap k_2^{\ast}=k, \quad
G_0=\Gal(k_1^{\ast}/k),\quad G_p=\Gal(k_2^{\ast}/k).
\end{gather}

Thus, $k_1^{\ast}/k$ is tamely ramified and $k_2^{\ast}/k$ 
is a $p$--extension so that it is wildly
ramified unless it is an extension of constants.

Now we study the field $k_1^{\ast}$. 
To find an explicit description of $k_1^{\ast}$ we proceed as
follows. Let
\begin{gather}\label{Eq3.1'}
F_0:=\prod_{\P\in{\ma P}_k^{\ast}}F_{\P}=\prod_{i=1}^t F_{\P_i}
\end{gather}
where $k\subseteq F_{\P}\subseteq \lam P$ is the unique subfield of 
the extension $\lam P/k$ of degree 
\begin{gather}\label{Eq3.3'}
c_{\P}:=\gcd (e_{\P},q^{d_{\P}}-1)=\gcd (e^{(0)}_{\P},q^{d_{\P}}-1).
\end{gather}
 Therefore $F_0$
satisfies that $KF_0/K$ is unramified at every finite
prime (Proposition \ref{P3.1}). 

Let $R:=k(\Lambda_{P_1\cdots P_t})$ and $R^+:=k(
\Lambda_{P_1\cdots P_t})^+$. Then $F_0\subseteq R$.

\begin{theorem}\label{T3.4'} With the notations as above, we have
\[
k_1^{\ast}\subseteq F_0{\ma F}_{q^{u_1}}\quad \text{and}\quad
K(F_0\cap R^+){\ma F}_{q^{t_0^{\prime}}} \subseteq K k_1^{\ast}
\subseteq KF_0 {\ma F}_{q^{u_1}},
\]
for some $u_1\in{\ma N}$ and
where $F_0$ is given by {\rm (\ref{Eq3.1'})}, $t_0$
is given by Proposition {\rm{\ref{P3.2}}} and $t_0=t_0^{\prime}
p^v$ with $\gcd(t_0^{\prime}, p)=1$.
\end{theorem}

\proof 
We will prove that $k_1^{\ast}\subseteq F_0{\ma F}_{q^{u_1}}$ for some
$u_1\in{\ma N}$.
For any prime $\P\in {\ma P}_k^{\ast}$
we obtain from Proposition \ref{P3.1} that if the ramification
index of $\P$ in $k_1^{\ast}/k$ is $b_{\P}$, then 
$b_{\P}|e_{\P}$ and since $k_1^{\ast}/k$ is a finite abelian
tamely ramified 
extension, we have $b_{\P}\mid q^{d_{\P}}-1$
(Proposition \ref{PalestineP4.1}). Hence
$b_{\P}\mid 
c_{\P}=\gcd(e_{\P},q^{d_{\P}}-1)=[F_{\P}:k]$. 
Let $F_{\P}^{\prime}$ be the
subfield of $F_{\P}$ of degree $b_{\P}$ over $k$.

We may assume that the finite ramified primes in $k_1^{\ast}/k$
are all of $\P_1,\ldots,
\P_t$ since, if some of the $b_{\P_i}$ are equal to $1$,
the argument below works even in this case.

We start with $\P_1$. From Abhyankar's Lemma we have that the
ramification index of $\P_1$ in $k_1^{\ast}F_{\P_1}^{\prime}$ over $k$ is 
$b_{\P_1}$. Let $I_{\P_1}$ be the inertia group of $\P_1$
in $k_1^{\ast} F^{\prime}_{\P_1}$ which is of order $b_{\P_1}$.
Let $E_1$ be the fixed field of $k_1^{\ast}F^{\prime}_{P_1}$ under
$I_{\P_1}$. Since $\P_1$ is fully ramified in $F^{\prime}_{\P_1}/k$
and unramified in $E_1/k$ we have $E_1\cap F^{\prime}_{\P_1}=k$
and 
\begin{gather*}
[E_1F^{\prime}_{\P_1}:k]=[E_1:k][F^{\prime}_{\P_1}:k]=
\frac{[k_1^{\ast}F^{\prime}_{\P_1}:k]}{|I_{\P_1}|}|I_{\P_1}|=
[k_1^{\ast}F^{\prime}_{\P_1}:k].
\end{gather*}
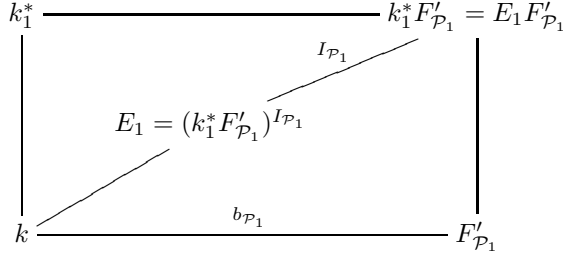
\[
\xymatrix{k_1^{\ast}\ar@{-}[rr]
\ar@{-}[dd]&&k_1^{\ast}F^{\prime}_{\P_1}=
E_1F^{\prime}_{\P_1}\ar@{-}[dd]
\ar@{-}[dl]_{I_{\P_1}}\\ &E_1=(k_1^{\ast}F^{\prime}_{\P_1})^{I_{\P_1}}
\ar@{-}[dl]\\
k\ar@{-}[rr]^{b_{\P_1}}&&F^{\prime}_{\P_1}}
\]

Therefore $k_1^{\ast}F^{\prime}_{\P_1}=E_1F^{\prime}_{\P_1}$.
Furthermore, since
$\P_2, \ldots, \P_t$ are unramified in $F^{\prime}_{\P_1}$ their
ramification indices are $b_{\P_2},\ldots, b_{\P_t}$ 
in $E_1F^{\prime}_{\P_1}/F^{\prime}_{\P_1}$.
Thus $\P_2, \ldots, \P_t$ 
have ramification indices $b_{\P_2},\ldots, b_{\P_t}$ 
in $E_1/k$.

Take now $E_1$ instead of $k_1^{\ast}$ and $F^{\prime}_{\P_2}$
instead of $F^{\prime}_{\P_1}$.
We obtain $E_2$ such that 
$E_1F_{\P_2}^{\prime}=E_2F_{\P_2}^{\prime}$ and $\P_3,
\ldots,\P_t$ are the only finite primes of $k$ ramified
in $E_2$ with ramification indices $b_{\P_3},\ldots b_{\P_t}$
respectively. Note that 
\[
k_1^{\ast}F^{\prime}_{\P_1}F^{\prime}_{\P_2}=
E_1F^{\prime}_{\P_1}F^{\prime}_{\P_2}=F^{\prime}_{\P_1}E_1
F^{\prime}_{\P_2}=F^{\prime}_{\P_1}E_2F^{\prime}_{\P_2}=
E_2F^{\prime}_{\P_1}F^{\prime}_{\P_2}.
\]

In the general
step we have $E_{i-1}F_{\P_i}^{\prime}=E_iF_{\P_i}^{\prime}$ and the
ramification indices of $\P_{i+1},\ldots, \P_t$ in $E_i/k$ are
 $b_{\P_{i+1}}, \ldots, b_{\P_t}$
and $k_1^{\ast}F^{\prime}_{\P_1}\ldots F^{\prime}_{\P_i}=
E_i F^{\prime}_{\P_1}\ldots F^{\prime}_{\P_i}$.

Keeping on in this way we finally obtain 
$E_t$ such that 
$E_{t-1}F_{\P_t}^{\prime}=E_tF_{\P_t}^{\prime}$,
no finite prime is ramified in $E_t/k$ and
$k_1^{\ast} F_0^{\prime} = E_t F_0^{\prime}$ where
$F_0^{\prime}=\prod_{i=1}^t F^{\prime}_{\P_i}$.

Since the only possibly ramified
prime in $E_t/k$ is $\p$ and it is tamely ramified,
from Proposition \ref{PalestineP4.2}
we obtain that $E_t/k$ is a constant field extension,
say $E_t={\ma F}_{q^{u_1}}(T)=k_{u_1}$.

Since $\big\{F_{\P_i}^{\prime}\big\}_{i=1}^t$ are pairwise
linearly disjoint and $F_0^{\prime}/k$ is a geometric extension, we have
\[
[F_0^{\prime}:k]=\prod_{i=1}^t [F_{\P_i}^{\prime}:k]=\prod_{i=1}^t
b_{\P_i}, \quad E_t\cap F_0^{\prime}=k \quad\text{and}\quad 
[k_1^{\ast}F_0^{\prime}:k]=[E_t:k][F_0^{\prime}:k].
\]
In particular, ${\ma F}_{q^{u_1}}$ is the field of constants of
$k_1^{\ast}F_0^{\prime}$.


Therefore $k_1^{\ast}\subseteq k_1^{\ast}F_0^{\prime}= E_t F_0^{\prime} \subseteq  F_0{\ma F}_{q^{u_1}}$ and
$Kk_1^{\ast}\subseteq KF_0{\ma F}_{q^{u_1}}$.
Finally, since the extension $K(F_0\cap R^+){\ma F}_{q^{t_0}}/K$
is unramified and the infinite primes are fully decomposed, it follows
that $K(F_0\cap R^+){\ma F}_{q^{t_0^{\prime}}}
\subseteq Kk_1^{\ast}$. \fin

\begin{remark}\label{R3.4''}{\rm{
In the proof of Theorem \ref{T3.4'} we have obtained that in fact
$k_1^{\ast}\subseteq E_tF_0^{\prime}$ and that ${\ma F}_{q^{u_1}}$ is the
field of constants of $k_1^{\ast}F_0^{\prime}$.
}}
\end{remark}

To study $k_2^{\ast}$ we first prove:

\begin{lemma}\label{L3.5}
We have $\g {k^{\ast}} =\g {({k^{\ast}_1})}
\g {({k^{\ast}_2})}= k^{\ast}$. Furthermore $\g {({k^{\ast}_1})}
=k^{\ast}_1$ and $\g {({k^{\ast}_2})}=k^{\ast}_2$.
\end{lemma}

\proof We have $k^{\ast}=k_1^{\ast}k_2^{\ast}$ and we have
already noted that $\g {k^{\ast}}=k^{\ast}$. 
Since $\g {({k^{\ast}_1})}/k_1^{\ast}$ is unramified and the
infinite primes decompose fully, the same holds in the 
extension $k^{\ast}
\g {({k^{\ast}_1})}/k^{\ast}$ so that $\g {({k^{\ast}_1})}\subseteq
\g {k^{\ast}}$. Similarly $\g {({k^{\ast}_2})}\subseteq \g {k^{\ast}}$.
Hence $\g {({k^{\ast}_1})}\g {({k^{\ast}_2})}\subseteq \g {k^{\ast}}$.

Now, since $\g {({k^{\ast}_1})}\supseteq k_1^{\ast}$ and 
$\g {({k^{\ast}_2})}\supseteq k_2^{\ast}$, we obtain
\[
\g {k^{\ast}}=k^{\ast} =k_1^{\ast} k_2^{\ast}\subseteq 
\g {({k^{\ast}_1})}\g {({k^{\ast}_2})}\subseteq \g {k^{\ast}}.
\]

Let now $[k_1^{\ast}:k]=a$ and $[k_2^{\ast}:k]=p^v$ where $p\nmid a$.
If $k_1^{\ast}\subsetneqq \g {({k^{\ast}_1})}$, let
$M:=\g {({k^{\ast}_1})}\cap k_2^{\ast}$. From the Galois correspondence
we obtain that $M\neq k$.

Let $[M:k]=p^b$ with $b\geq 1$. We have $M/k$ is unramified since
otherwise there exists a prime in $k$ with ramification index $p^{c}$ with
$c\geq 1$ in $M$. Since $p\nmid a$, it follows that there exists a 
ramified prime in $\g {({k^{\ast}_1})}/k_1^{\ast}$ with ramification index
$p^c$. This contradiction shows that $M/k$ is unramified. Thus $M/k$
is an extension of constants. It follows that $\p$ has inertia degree
$p^b$ in $M/k$ but this implies that the inertia degree of the
infinite primes in $\g {({k^{\ast}_1})}/k_1^{\ast}$ is $p^b$
which is impossible.
Therefore $\g {({k^{\ast}_1})}=k_1^{\ast}$.
Similarly $\g {({k^{\ast}_2})}=k^{\ast}_2$. \fin

\begin{remark}\label{R3.6}{\rm{
In general, if $L=L_1L_2$, then $\g {(L_1)}\g {(L_2)} \subseteqq
\g L$ but not necessarily $\g L=
\g {(L_1)}\g {(L_2)}$. For instance, let $q>2$ and
let $P,Q,R,S\in R_T$ be four monic
irreducible polynomials in $k$. Let $L_1:=k(\Lambda_{P^2Q^2})^+$
and $L_2:=k(\Lambda_{R^2S^2})^+$. Then $L_1=\g {(L_1)}$ and
$L_2=\g {(L_2)}$. Let $L:=L_1L_2$. Then $\g L=k(\Lambda_{
P^2Q^2R^2S^2})^+$ and $[\g L:L]=q-1>1$. Thus $\g L=\g {(L_1
L_2)}\neq \g {(L_1)}\g {(L_2)}= L$.
}}
\end{remark}

Now let $\Gal(k_2^{\ast}/k)\cong C_{p^{n_1}}\times \cdots
\times C_{p^{n_{\nu}}}$ and if for each $1\leq i\leq \nu$, $E_i$
is the subfield $k\subseteq E_i\subseteq k_2^{\ast}$
such that $\Gal(E_i/k)\cong C_{p^{n_i}}$, then from
\cite[Theorem 5.7]{MaRzVi2013}, we obtain that $\g {(E_i)}$ is the
composite of $p$--cyclic extensions of $k$ such that in each
one only one prime is ramified or is an extension of constants.
Therefore, $\g {({k^{\ast}_2})}=k_2^{\ast}$ is the composite
of this type of cyclic $p$--extensions.

Finally, since $u_{\P_j}\geq 1$ for $s+1\leq j\leq t$ (see (\ref{Eq0.1})),
we have the following result.

\begin{theorem}\label{T3.6}
The field $k_2^{\ast}$ is of the form 
$k_2^{\ast} = J_{s+1} J_{s+2}\cdots J_t J_{\infty}$ where
$\P_j$ is the only ramified prime in $J_j/k$, 
$[J_j:k]=p^{v_j}$ with $0\leq v_j\leq u_{\P_j}$ for $s+1\leq j\leq t$,
and $J_{\infty}$ is an abelian $p$--extension of $k$ that is
either an extension of constants or 
such that $\p$ is the only ramified prime. \fin
\end{theorem}

\section{The genus field in a special case}\label{S3}

Let $K/k$ be a finite separable extension
such that for all $\P\in{\ma P}_k$, $p\nmid e_{\P}=
\gcd(e_1,\ldots,e_r)$ where $\con_{k/K}\P=\producto$. 
That is, we assume $t=s$ and also $p\nmid e_{\p}$.

We have that $k_1^{\ast}$ is given by (\ref{Eq3.1''}) and
in this case we have $k_2^{\ast}/k$ is unramified.
Hence $k_2^{\ast}/k$ is
an extension of constants. 

To find a more explicit description of $k_1^{\ast}$ we proceed as
follows. First we consider the behavior of $\p$.
Let $\con_{k/K}\p$ be given by (\ref{Eq3.2'}).

We have
$e_{\infty}(F_{\P}|k)\mid
\gcd(c_{\P},q-1)$ for $\P\in {\ma P}_k^{\ast}$. By
Abhyankar's Lemma we obtain that if 
\[
c_{\infty}:=e_{\infty}(F_0|k),
\]
is the ramification index of $\p$ in $F_0/k$
then
\[
c_{\infty}\mid \lcm\big[\gcd(e_{\P_1},q-1),\ldots,\gcd(e_{\P_s},q-1)\big]=
\gcd\big(\lcm[e_{\P_1},\ldots, e_{\P_s}],q-1\big).
\]

To obtain a formula for $c_{\infty}$, we consider the following.
We have from (\ref{Eq3.3'})
\[
c_{\P}=[F_{\P}:k]=\gcd (e_{\P},q^{d_{\P}}-1),
\]
where $d_{\P}=d_k(\P)$. Let $H:=\Gal(R/F_0)$, where $R=k(
\Lambda_{P_1\cdots P_s})$. Let $M:=F_0R^+$, where $R^+=k(
\Lambda_{P_1\cdots P_s})^+$.
Therefore
\begin{align*}
e_{\infty}(M|k)&=e_{\infty}(M|R^+)e_{\infty}(R^+|F_0\cap R^+)
e_{\infty}(F_0\cap R^+|k)\\
&=[M:R^+]\cdot 1\cdot 1=[M:R^+]=
[F_0:F_0\cap R^+];\\
e_{\infty}(M|k)&=e_{\infty}(M|F_0)e_{\infty}(F_0|F_0\cap R^+)
e_{\infty}(F_0\cap R^+|k)\\
&=1\cdot e_{\infty}(F_0|F_0\cap R^+)\cdot 1=
e_{\infty}(F_0|F_0\cap R^+).
\end{align*}
Hence
\begin{gather}\label{Eq3.2}
c_{\infty}=e_{\infty}(F_0|k)=e_{\infty}(F_0|F_0\cap R^+)=e_{\infty}
(M|k) =[F_0:F_0\cap R^+]=[M:R^+].
\end{gather}

We choose the maximal field $F$ with $F_0\cap R^+\subseteq F
\subseteq F_0$ and such that the
infinite primes of $K$ decompose fully in $KF$.
Note that such field $F$ exists since if $F_1,F_2$ are two fields such
that $F_0\cap R^+\subseteq F_i\subseteq F_0$ and such that the
infinite primes of $K$ decompose fully in $KF_i/K$, $i=1,2$, then
$F_1F_2$ satisfies the same properties.

\begin{remark}\label{R4.-1}{\rm{
With the notation
of Theorem \ref{T3.4'},
observe that since $KF/K$ is unramified, and
because $\p$ splits fully in $KF/K$,
it follows that $F\subseteq k_1^{\ast}$  
so that $Fk_1^{\ast}=k_1^{\ast}
\subseteq k_1^{\ast}F_0^{\prime}$. Since $F_0\cap
R^+\subseteq F_0^{\prime}\subseteq F_0$ we have 
$F\subseteq F_0^{\prime}$ .
In general we may have $F_0^{\prime}\neq F$, 
see Example \ref{Ex5.1}.
}}
\end{remark}

Next, we determine $F$ for an abelian extension $K/k$.

\begin{proposition}\label{P4.0} 
Let $K/k$ be a finite abelian tamely ramified extension. With the
notation in Theorem {\rm{\ref{T2.6}}} we have
\begin{gather*}
F\subseteq \g E \subseteq F_0,
\intertext{more precisely}
F=\g {E^{H_1}}\quad\text{and}\quad \g K=KF.
\end{gather*}
\end{proposition}

\proof In this case $s=t$ and $N=P_1\cdots P_t$. Since
for any prime $\P$ in $k$, the ramification index in $K/k$ is the
same as the ramification index in $E/k$ (see \cite[Section 4.1]
{MaRzVi2013}) and $F_0=\prod_{\P\in{\ma P}_k^{\ast}} F_{\P}$,
we have $\g E\subseteq F_0$.

The infinite prime decomposes fully in $F_0\cap R^+/k$. Hence
the infinite primes decompose fully in $E(F_0\cap R^+)/E$. Since
the extension $E(F_0\cap R^+)/E$ is unramified, we have
$F_0\cap R^+\subseteq \g E$. 

We observe that by Abhyankar's Lemma (Theorem \ref{T2.3}) the
extension $K(F_0\cap R^+)/K$ is unramified and the infinite primes
decompose fully. Thus $F\subseteq \g E$.

Finally, again by Abhyankar's Lemma, $K\g E/K$ is unramified and the
inertia of the infinite primes corresponds to $H_1$, that is,
$\g {E^{H_1}}$ is the maximal
extension such that $F_0 \cap R^+\subseteq \g {E^{H_1}}\subseteq
F_0$ and that in $K\g {E^{H_1}}/K$ the infinite primes decompose
fully. Therefore $F=\g {E^{H_1}}$. From Theorem \ref{T2.6},
it follows $\g K=KF$. $\fin$

\begin{remark}\label{Ex5.4}{\rm{
Let $K/k$ be a finite tamely ramified abelian extension. Let $\P_1,
\ldots,\P_s$ be the finite ramified primes. Then
$F_0=\prod_{i=1}^s F_{\P_i}$ with $k\subseteq F_{\P_i}
\subseteq k(\Lambda_{P_i})$. We have $[F_{\P_i}:k]=c_{\P_i}=\gcd
(e_{\P_i},q^{\deg P_i}-1)$. Since $K/k$ is abelian and tamely
ramified we have $e_{\P_i}\mid q^{\deg P_i}-1$ 
(Proposition \ref{PalestineP4.1}). Therefore, $c_{\P_i} = e_{\P_i}$. 
}}
\end{remark}

Now let 
\[
c^{\prime}_{\infty}:=[F:F_0\cap R^+]=e_{\infty}(F|k).
\]
Since 
\[
c_{\infty}=[F_0:F_0\cap R^+] =[F_0:F][F:F_0\cap R^+]=
[F_0:F]c^{\prime}_{\infty}
\]
we have $c^{\prime}_{\infty}\mid c_{\infty}$.

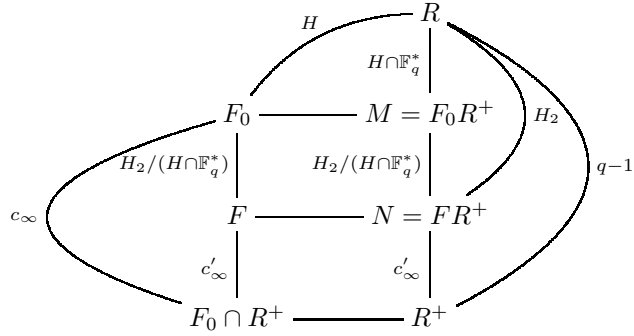
\[
\xymatrix{&R\ar@{-}[d]_{H\cap \f}
\ar@/^5pc/@{-}[ddd]^{q-1}\ar@/^3pc/@{-}[dd]^{H_2}\\
F_0\ar@{-}[r]\ar@{-}[d]_{H_2/(H\cap \f)}\ar@/^1pc/@{-}[ur]^H
\ar@/_6pc/@{-}[dd]_{c_{\infty}}&
M=F_0 R^+\ar@{-}[d]_{H_2/(H\cap \f)}\\
F\ar@{-}[r]\ar@{-}[d]_{c^{\prime}_{\infty}}&N=FR^+\ar@{-}[d]_{
c^{\prime}_{\infty}}\\
F_0\cap R^+
\ar@{-}[r]&R^+}
\]

The degree $c^{\prime}_{\infty}$ must satisfy the following. By
Abhyankar's Lemma, we have that if ${\eu P}$ is a prime in $KF$
dividing $\p$ and if ${\eu P}\cap K={\eu p}_{i,\infty}$, then
\begin{align*}
e({\eu P}|\p)&=\lcm[e_{i,\infty},c^{\prime}_{\infty}]=\frac{
e_{i,\infty}c^{\prime}_{\infty}}{\gcd(e_{i,\infty},c^{\prime}_{\infty})}\\
&=e({\eu P}|{\eu p}_{i,\infty})e({\eu p}_{i,\infty}|\p)=
e({\eu P}|{\eu p}_{i,\infty})e_{i,\infty}.
\end{align*}
It follows that 
\begin{gather}\label{Eq3.3}
e({\eu P}|{\eu p_{i,\infty}})=\frac{c^{\prime}_{\infty}}
{\gcd(e_{i,\infty},c^{\prime}_{\infty})}.
\end{gather}
Therefore
\begin{gather*}
e({\eu P}|{\eu p}_{i,\infty})=1\iff
\gcd (e_{i,\infty},c^{\prime}_{\infty})=c^{\prime}_{\infty} 
\iff c^{\prime}_{\infty}\mid e_{i,\infty}.
\end{gather*}
Thus, $KF/K$ is unramified if and only if
$c^{\prime}_{\infty}\mid e_{\P_{\infty}}=\gcd (e_{1,\infty},\ldots,
e_{r_{\infty},\infty})$.

Therefore $c^{\prime}_{\infty}$ must be maximal
in the sense $c^{\prime}_{\infty}\mid c_{\infty}$,
$c^{\prime}_{\infty}\mid 
e_{\infty}$ where $e_{\infty}=e_{\p}$,  
$c_{\infty}$ is given by (\ref{Eq3.2})
and the infinite primes of $K$ decompose
fully in $KF$. Hence 
\begin{gather}\label{Eq3.30}
c^{\prime}_{\infty}\mid\gcd (c_{\infty}, e_{\infty}).
\end{gather}
That is, $F$ is the field
\begin{gather}\label{Eq3.31}
F_0\cap R^+\subseteq F\subseteq F_0\quad\text{such that}\quad
[F:F_0\cap R^+]=c^{\prime}_{\infty}.
\end{gather}

Let $H_2$ be the subgroup of $\f$ of order $\frac{q-1}
{c_{\infty}^{\prime}}$ and let $N:=R^{H_2}$. Note that 
$|H\cap \f|=[R:F_0R^+]=
\frac{q-1}{c_{\infty}}$. Hence $|H\cap \f|\mid |H_2|$ and from
(\ref{Eq3.30}) we obtain
\[
[M:N]=[F_0:F]=\frac{c_{\infty}}{c_{\infty}^{\prime}}.
\]

With the above notation we have the following result.

\begin{theorem}\label{T3.4} Let $K/k$ be a finite separable 
extension such that
every prime $\P\in {\ma P}_k$ satisfies that if
$\con_{k/K}\P=\producto$, then $p\nmid e_{\P}= \gcd(e_1,\ldots,e_r)$.
Then 
\[
KF{\ma F}_{q^{t_0}} \subseteq \g K \subseteq KF_0 {\ma F}_{q^{u}},
\]
where $F_0$ is given by {\rm (\ref{Eq3.1'})}, F is 
given by {\rm (\ref{Eq3.31})}, $t_0$ 
is given by {\rm{(\ref{Eq3.2''})}}
and $u\in{\ma N}$.

Furthermore, $KF_0{\ma F}_{q^{u}}/K$ is unramified
at every finite prime and the ramification index of the 
infinite prime ${\eu p}_{i,\infty}$
is $\frac{c_{\infty}}{\gcd(e_{i,\infty},c_{\infty})}$, $1\leq i\leq r_{\infty}$
where $c_{\infty}$ is given by {\rm{(\ref{Eq3.2})}}.

\end{theorem}

\proof 
Because $KF/K$ is unramified and the infinite
primes decompose fully, we have $F\subseteq k_1^{\ast}$.
Therefore by Proposition \ref{P3.2} we have
 $KF\finite {t_0}\subseteq K k^{\ast}=\g K$.

Since $p\nmid e_{\P}$ for all $\P\in{\ma P}_k$ it follows
from Theorem \ref{T3.6} and 
Proposition \ref{PalestineP4.2} that $k_2^{\ast}/k$
is an extension of constants, so that $k_2^{\ast}={\ma F}_{
q^{u_2}}(T)$. Furthermore, since $k_2^* \subseteq \g K$ we have $u_2|t_0$. 

From Theorem \ref{T3.4'} we have $k_1^{\ast}\subseteq F_0
{\ma F}_{q^{u_1}}$ for some $u_1\in{\ma N}$
and $Kk_1^{\ast}\subseteq KF_0{\ma F}_{q^{u_1}}$.
Hence $\g K=Kk^{\ast}=Kk_1^{\ast}k_2^{\ast}\subseteq
KF_0{\ma F}_{q^{u}}$, where $u=\lcm [u_1,u_2]$.

The ramification index of $\p$ in
$F_0/k$ is $c_{\infty}$ where $c_{\infty}$ is
given by (\ref{Eq3.2}). 
Now applying (\ref{Eq3.3}) to $F_0$ and $c_{\infty}$
we obtain the ramification index of ${\eu p}_{i,\infty}$
in $KF_0{\ma F}_{q^{u}}/K$ is $\frac{c_{\infty}}
{\gcd(e_{i,\infty},c_{\infty})}$.
$\fin$

\begin{remark}\label{R3.5} {\rm{
Note that 
$F\subseteq \g K\cap F_0$. 
Because $\g K\cap F_0 \subseteq \g K$, we have the infinite primes of $K$ decompose
fully in $(\g K\cap F_0 )K$. Besides, $F_0 \cap R^+ \subseteq \g K\cap F_0 \subseteq F_0$. It follows from the maximality
of $F$ that $$\g K\cap F_0=F.$$ Therefore $KF\finite {t_0}
\cap F_0=F$ and $\g K\cap KF_0\finite {t_0}=KF\finite {t_0}$.

Observe that if we had in the proof that $
(\g K)_u\cap F_0=F$, then by the Galois correspondence we
would obtain $(\g K)_u=((\g K)_u\cap F_0) K_u=FK{\ma F}_{q^u}$.
Hence $FK{\ma F}_{q^{t_0}} \subseteq \g K\subseteq (\g K)_u=
FK{\ma F}_{q^u}=KF{\ma F}_{q^{t_0}}{\ma F}_{q^u}$.
Therefore $\g K/FK{\ma F}_{q^{t_0}}$ is an extension of
constants and since the field of constants of $\g K$ is
${\ma F}_{q^{t_0}}$, it would follow the equality
\[
\g K =FK {\ma F}_{q^{t_0}}.
\]

In case that $F=F_0$ then $KF{\ma F}_{q^{t_0}}\subseteq
\g K\subseteq KF{\ma F}_{q^u}$ so that $\g K/KF{\ma F}_{q^{t_0}}$
is an extension of constants and then $\g K= KF{\ma F}_{q^{t_0}}$.

Finally, if $u=t_0$, then $\g K=\g K \cap KF_0{\ma F}_{q^{t_0}}=
KF{\ma F}_{q^{t_0}}$ (see the diagram below). Hence 
$\g K=KF{\ma F}_{q^{t_0}}$.

Also, when $K/k$ is an abelian tamely ramified extension, we
have $\g K= KF$ (see 
Proposition \ref{P4.0}).

In short, it is very likely that always $\g K=KF{\ma F}_{q^{t_0}}$.

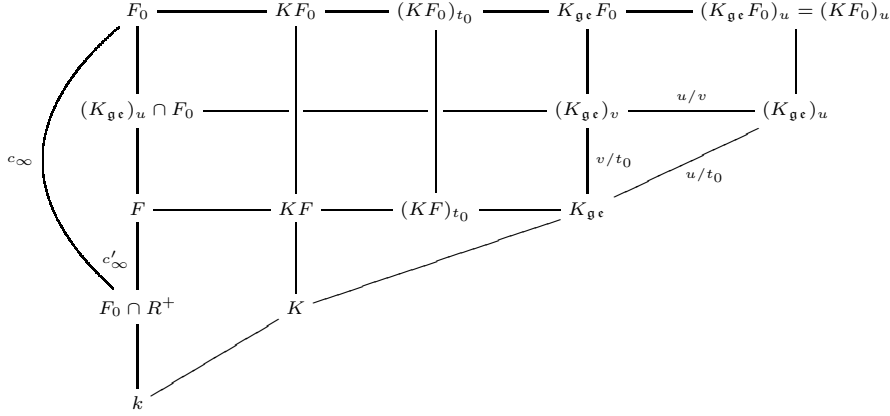
\begin{scriptsize}
\[
\xymatrix{F_0\ar@{-}[d]\ar@{-}[r]\ar@/_3pc/@{-}[ddd]_{c_{\infty}}
&KF_0\ar@{-}[dd]\ar@{-}[r]&(KF_0)_ {t_0}\ar@{-}[r]
\ar@{-}[dd]&\g K F_0\ar@{-}[d]\ar@{-}[r]&
(\g K F_0)_u=(KF_0)_u\ar@{-}[d]\\
({\g K})_u\cap F_0\ar@{-}[d]\ar@{-}[rrr]
|!{[r];[dd]}\hole|!{[rr];[dd]}\hole&&& ({\g K})_v\ar@{-}[r]^{u/v}
\ar@{-}[d]^{v/t_0}&({\g K})_{u}\\
F\ar@{-}[r]\ar@{-}[d]_{c^{\prime}_{\infty}}
&KF\ar@{-}[r]\ar@{-}[d]&(KF)_{t_0}\ar@{-}[r]
&\g K\ar@{-}[ur]_{u/t_0}\ar@{-}[dll]\\
F_0\cap R^+\ar@{-}[d]&K\ar@{-}[dl]\\
k}
\]
\end{scriptsize}
}}
\end{remark}

\section{Applications and examples}\label{S5}

\begin{example}\label{Ex5.1} {\rm{Consider 
$q=3$ and $P=T^3+2T+1$. We have
$P$ is irreducible in ${\ma F}_3(T)$. Let $K=k(\sqrt{P})$. In our
construction, if $\P$ is the prime corresponding to $P$, we have
$F_0=F_{\P}=k(\sqrt{(-1)^{\deg P}P})=k(\sqrt{-P})$. Now
$\p$ is ramified in $K$ and in $F_0=F_{\P}$. Therefore 
$t_0=1$, that is, the field of constants of $\g K$ is ${\ma F}_3$.
Since $[R^+:k]=13$ and $[F_0:k]=2$, we have $F_0\cap R^+
=k$. Now $KF_0=K(\sqrt{-1})$. Since $\sqrt{-1}\notin
{\ma F}_3$ we have $K(\sqrt{-1})=K{\ma F}_9$ and the infinite
primes are inert in $KF_0/K$. Hence $F=k$ and $\g K=K$.
Here we have $F_0^{\prime}=F_0=k(\sqrt{-P})\neq k=F$.
}}
\end{example}

\subsection{Cyclic extensions of prime degree not dividing $q(q-1)$}\label{S5.1}

Let $l$ be a prime not dividing $q(q-1)$ and let $K/k$ be a cyclic extension
of degree $l$. Let $\P_1,\ldots, \P_t$ be the primes in $k$
ramified in $K$. Note that since the inertia group of a ramified prime
is contained in the multiplicative group of the residue field, we have
$l\mid (q^{\deg P_i}-1)$ for $1\leq i\leq t$. In particular
$\p$ is not ramified. In this case we have
$k\subseteq F_{\P i}\subseteq \lam {{P_i}}$
for $1\leq i\leq t$ where $F_{\P i}$ is
the unique subfield of $\lam {{P_i}}$ of degree $c_{\P_i}=\gcd(
e_{\P_i},q^{d_{\P_i}}-1)=l$. Then
\[
F_0=\prod_{i=1}^t F_{\P i}\subseteq k(\Lambda_{P_1\cdots P_t})^+.
\]
Therefore we have $c_{\infty}=e_{\p}=1$, $F_0\cap R^+=F_0$.
Thus $F=F_0$, $c_{\infty}^{\prime}=1$. Furthermore, if $t_0$
is the degree of the infinite prime(s) above $\p$ in $K$, then
$t_0=1$ or $l$. In fact $t_0=1$ iff $\p$ decomposes in $K/k$.
This is equivalent to
$K\subseteq k(\Lambda_{P_1\cdots P_t})^+$. We have $t_0=l$
iff $\p$ is inert in $K/k$ iff
$K\nsubseteq k(\Lambda_{P_1\cdots P_t})^+$.

From Proposition \ref{P3.2} and Theorem \ref{T3.4}
we have $\g K {}=KF{\ma F}_{q^{t_0}}$, since in this case
we have $F=F_0$ and $u=t_0$. 

First we consider
$K\subseteq k(\Lambda_{P_1\cdots P_t})^+$. Then $\g K {}=
KF{\ma F}_{q^{t_0}}=KF=F$ and $[\g K {}:K]=l^{t-1}$.

Now we consider $K\nsubseteq k(\Lambda_{P_1\cdots P_t})^+$.
Then $K\nsubseteq F$ and in particular $k\subseteq K\cap F
\subsetneqq K$ so that $K\cap F=k$ and $[KF:K]=[F:k]=l^t$. 
We will prove ${\ma F}_{q^l}\subseteq KF$. First, we have
$[KF:k]=[{\ma F}_{q^l} F:k]=l^{t+1}$. Now if $k_l:={\ma F}_{q^l}(T)$,
then $k_l
\cap K=k$. Now $\p$ is inert in $K/k$ and in $k_l/k$. The
decomposition group ${\mathcal D}$ of $\p$ in $K_l=Kk_l$ is 
a cyclic group of order
$l$. Consider $L:=(K_l)^{\mathcal D}$. The prime $\p$ decomposes fully
in $L/k$ and $\P_1,\ldots,\P_t$ are the ramified primes in $L/k$.
It follows that $L\subseteq F$. Since $L\neq K$ we obtain that
$KL=K_l$ and $KL=K_l=K{\ma F}_{q^l}\subseteq KF$. Thus
${\ma F}_{q^l}\subseteq KF$.
Therefore $\g K {}=KF{\ma F}_{q^l}=KF$ and $[\g K {}:K]=
[KF:K]=l^t$.

Note that this example is consistent with the results of 
\cite{MaRzVi2013, MaRzVi2015}, in particular with 
Theorem 4.2 and Remark 4.3 of \cite{MaRzVi2015}. In the notation
of those papers, $K\subseteq F\iff E=K$ and $[\g K {}:K]=
[\g E {}:E]=l^{t-1}$. When $K\nsubseteq F$, $\g E{}=F$
we have $[\g K{}:K]=l^t=l l^{t-1}=l[\g E {}:E]$.

\subsection{Radical extensions}\label{S5.2}

Let $K=k(\sqrt[n]{\gamma D})$, where $D\in R_T$ is a monic polynomial
and $\gamma \in
\f$. Let $D=P_1^{\alpha_1}
\cdots P_s^{\alpha_s}$ be the decomposition of $D$ as product of
irreducible polynomials. We assume that $D$ is $n$--th power free, that is, 
$0<\alpha_i<n$ for $1\leq i\leq s$ and we also assume $p\nmid n$.

The finite ramified primes are $\P_1,\ldots,\P_s$ and they are tamely ramified.
Indeed, let ${\eu p}_i$ be a prime in $K$ above $\P_i$. We have
\begin{gather}\label{Eq5.1}
e({\eu p}_i|\P_i) v_{\P_i}(D)=e({\eu p}_i|\P_i) \alpha_i=
v_{{\eu p}_i}(D)=v_{{\eu p}_i}((\sqrt[n]{\gamma D})^n)=nv_{{\eu p}_i}
(\sqrt[n]{\gamma D}).
\end{gather}
Let $d_i=\gcd(\alpha_i,n)$. We obtain from (\ref{Eq5.1}) that
$\frac{n}{d_i}\mid e({\eu p}_i|\P_i)$. On the other hand if we write $K=k(y)$,
where $y^n=\gamma D$, set $z=y^{n/d_i}$. Thus $z^{d_i}=y^n =
\gamma D= \gamma (P_i^{\alpha_i/d_i})^{d_i} (D/P_i^{\alpha_i})$. Therefore
$k(z)=k(\sqrt[d_i]{\gamma D/P_i^{\alpha_i}})$. In particular $\P_i$
is unramified in $k(z)/k$. It follows that $e({\eu p}_i|\P_i)=n/d_i$.
Thus $e_{\P_i}=n/d_i$.
Similarly we obtain $e_{\infty}=e_{\infty}(K|k) = 
\frac{n}{d}$, where $d=\gcd(\deg D,n)$. 

Therefore $F_0=\prod_{i=1}^s F_{\P_i}$ with 
$c_{\P_i}=\gcd(e_{\P_i},q^{\deg P_i}-1)=
\gcd\big(\frac{n}{d_i},q^{\deg P_i}-1\big)$. Then 
\begin{gather*}
e_{\infty}(F_{\P_i}|k)\mid \gcd(c_{\P_i},q-1)=\gcd(e_{\P_i},q-1)
=\gcd\big(\frac{n}{d_i},q-1\big).\\
\intertext{Hence}
c_{\infty}\mid \gcd(\lcm[e_{\P_i},\ldots,e_{\P_s}],q-1)=\gcd(\lcm\big[
\frac{n}{d_1},
\ldots,\frac{n}{d_s}\big],q-1)=\gcd(\frac{n}{d_0},q-1),\\
\intertext{where $d_0=\gcd[d_1,\ldots,d_s]$. We also have}
c^{\prime}_{\infty}\mid \gcd(c_{\infty},e_{\infty})|\gcd\big(\frac{n}{d_0},
\frac{n}{d},q-1\big).
\end{gather*}

By Theorem \ref{T3.4} we obtain 
\[
KF{\ma F}_{q^{t_0}} =
k(\sqrt[n]{\gamma D})F{\ma F}_{q^{t_0}}\subseteq \g{K} =
\g {k(\sqrt[n]{\gamma D})}\subseteq KF_0{\ma F}_{q^u},
\]
where $t_0, u\in{\ma N}$.

To find $t_0$, 
we consider the subfield $E=k(\sqrt[d]{\gamma
D})\subseteq k(\sqrt[n]{\gamma D})$. Since $\p$ is unramified
in $E/k$ and fully ramified in $K/E$, we have that the 
inertia degree of $\p$ in $K/k$ is equal to the inertia degree
of $\p$ in $E/k$. Note that $D(T)=T^l+a_{l-1}T^{l-1}+\cdots+a_1
T+a_0=T^l\big(1+a_{l-1}(\frac{1}{T})+\cdots +a_1(\frac{1}{T})^{l-1}
+a_0 (\frac{1}{T})^l\big)=T^l D_1(\frac{1}{T})$ with $D_1(0)=1$
and $d|l$. Hence $E=k(\sqrt[d]{\gamma D_1(1/T)})$
with $D_1(1/T)\in {\ma F}_q[1/T]$ and $D_1(1/T)
\equiv  1 \bmod  (1/T) $. Therefore $X^d-\gamma D_1(1/T) \bmod\ \p$
becomes $\bar{X}^d-\gamma\in {\ma F}_q [\bar{X}]$.

Let $\mu\in\bar{\ma F}_q$ be a fixed $d$-th root of $\gamma$. If
$\zeta_d$ denotes a primitive $d$--th root of unity, we have that
the factorization of $\bar{X}^d-\gamma$ in ${\ma F}_q[\bar{X}]$ is
of the form
\[
\bar{X}^d-\gamma=\prod_{j=1}^r\Irr(\zeta_d^{i_j}\mu,\bar{X},{\ma F}_q)
\]
for some $0\leq i_1<i_2<\cdots <i_r\leq d-1$. From Hensel's Lemma,
we obtain $X^d-\gamma D_1\big(\frac{1}{T}\big)=\prod_{j=1}^r
F_j(X)$ with $F_j(X)\in k_{\infty}[X]$ distinct irreducible polynomials. 
In particular $\con_{k/K}\p={\eu p}_{\infty, 1}\cdots {\eu p}_{\infty,r}$
with $\deg_K {\eu p}_{\infty, j}=\deg F_j(X)$, $1\leq j\leq r$. Therefore
$t_0=\gcd_{1\leq j\leq r}\{\deg F_j(X)\}=\gcd_{1\leq j\leq r}
\{[{\ma F}_q(\zeta_d^{i_j}\mu):{\ma F}_q]\}=\gcd_{0\leq i\leq d-1}
\{[{\ma F}_q(\zeta_d^i \mu):{\ma F}_q]\}$. In short if we write $\sqrt[d]{
\gamma}=\mu$,
\begin{gather}\label{Eq5.1*}
t_0=\gcd_{0\leq j\leq r}\big[{\ma F}_q(\zeta_d^{i_j}\sqrt[d]{\gamma}):{\ma F}_q\big]
=\gcd_{0\leq i\leq d-1}\big[{\ma F}_q(\zeta_d^{i}\sqrt[d]{\gamma}):{\ma F}_q\big].
\end{gather}

It follows from (\ref{Eq5.1}) that if $\gcd(\alpha_i,n)=1$ 
for some $i$, then $K/k$ is a geometric extension.

\begin{example}\label{Ex5.3} {\rm{Consider 
$q=3$, $P_1 = T$, $P_2=T^2 - T -1$ and $D=P_1^2P_2$. We have
$P_1$ and $P_2$ are irreducible in ${\ma F}_3(T)$. Let $K=k(\sqrt[10]{-D})$. In our
construction, if $\P_i$ is the prime corresponding to $P_i$, we have 
$F_{\P_1}=k$ and $F_{\P_2}=k(\sqrt{P_2})$. Hence $F_0=k(\sqrt{P_2})$. 
On the other hand, $\p$ decomposes in $k(\sqrt{P_2})/k$, thus 
 $k(\sqrt{P_2}) \subseteq k(\Lambda_{P1P2})^+$. Therefore $F =F_0=k(\sqrt{P_2})$. 

Since $d=2$, we obtain $t_0 = 2$. Because $u_2$ is a power of $3$ and $u_2$ 
divides $t_0$, we have $u_2 = 1$. From the proof of Theorem \ref{T3.4'}, we obtain that in this case
$E_1 = k_1^*$ and ${\ma F}_{q^{u_1}} = E_2 \subseteq k_1^*$. Therefore $u_1$ divides $t_0$. Thus
$u = u_1 \in \{1,2 \}$. It follows from Theorem  \ref{T3.4} that $\g K = Kk(\sqrt{P_2}) {\ma F}_9$.
}}
\end{example}

\subsection{Radical extensions of prime power degree dividing
$q-1$}\label{S5.3}

As a particular case of
Subsection \ref{S5.2}, let $l$ be a prime number
such that $l^n\mid q-1$. Let $D\in R_T$ be a monic polynomial $l^n$--power
free. Let $D=P_1^{\alpha_1}\cdots P_s^{\alpha_s}$ with $P_1,\ldots,
P_s\in R_T^+$ and $v_l(\alpha_i)=a_i<n$. Let $\gamma \in \f$ and
$K=k(\sqrt[l^n]{\gamma D})$. Then $e_{\P_i}=l^{n-a_i}$, $1\leq i\leq s$.
Since $K/k$ is a cyclic extension of degree $l^n$, $K/k$ is a geometric
extension if and only if $a_i=0$ for some $1\leq i\leq s$.

Now, we have $F_{\P_i}\subseteq k(\Lambda_{P_i})$ and $c_{\P_i}=
\gcd(e_{\P_i},q^{\deg P_i}-1)= e_{\P_i}=l^{n-a_i}$. Therefore
$F_{\P_i}=k(\sqrt[l^{n-a_i}]{(-1)^{\deg P_i}P_i})$ and 
$F_0=\prod_{i=1}^s F_{\P_i}$.

We have $e_{\p}=e_{\infty}=l^{n-d}$, where $d=\min\{n,d^{\prime}\}$
and $v_l(\deg D)=d^{\prime}$. Furthermore the inertia degree of
$\p$ is $f_{\infty}=l^m$, where ${\ma F}_{q^{l^m}}=\F(\sqrt[l^d]
{(-1)^{\deg D}\gamma})$ (see \cite[Proposition 2.8]{BaRzVi2013}).
Hence $t_0=l^m$ and the field of constants of $\g K$ is ${\ma F}_{q^{l^m}}$.

Now, with respect to $\p$ we have 
$e_{\infty}(F_{\P_i}|k)=l^{n-a_i-d_i}$, where
$d_i=\min\{n-a_i,d_i^{\prime}\}$, $v_l(\deg P_i)=d^{\prime}_i$.
From Abhyankar's Lemma we obtain 
\[
e_{\infty}(F_0|k)=\lcm[e_{\infty}(F_{\P_i}|k)\mid 1\leq i\leq s]=
\lcm[l^{n-a_i-d_i}\mid 1\leq i\leq s]=l^{n-\delta},
\]
where $\delta=\min\limits_{1\leq i\leq s}\{a_i+d_i\}=
\min\limits_{1\leq i\leq s}\{a_i+\min\{n-a_i,d_i^{\prime}\}\}=
\min\limits_{1\leq i\leq s}\{n,v_l(\deg P_i^{\alpha_i})\}$. That is
\[
c_{\infty}=[F_0:F_0\cap R^+]=e_{\infty}(F_0|k)=l^{n-\delta}.
\]
Note that $d\geq \delta$. Then $c^{\prime}_{\infty}\mid\gcd(c_{\infty},
e_{\infty})=\gcd(l^{n-\delta},l^{n-d})=l^{n-d}$.
We have that $F$ is the subfield $F_0\cap R^+\subseteq F
\subseteq F_0$ such that $[F:F_0\cap R^+]=c_{\infty}^{\prime}\mid
l^{n-d}$.
\[
\xymatrix{
F_0\ar@{-}[d]^{\frac{c_{\infty}}{c_{\infty}^{\prime}}}
\ar@/_2pc/@{-}[dd]_{l^{n-\delta}=c_{\infty}}\\
F\ar@{-}[d]^{c^{\prime}_{\infty}\mid l^{n-d}}\\ F_0\cap R^+}
\]

\begin{example}\label{Ex5.3'}{\rm{
Let $k={\ma F}_5(T)$ and 
$K:=k\big(\sqrt[3]{T(T^2+T+1)}\big)\cdot {\ma F}_{5^2}=
k\big(T,\sqrt[3]{D(T)}\big)\cdot {\ma F}_{5^2}$, where $D(T)=T(T^2+T+1)$. Let
$K_0:=k(\sqrt[3]{T(T^2+T+1)})$. Note that $K=K_0(\zeta_3)
=K_0\cdot {\ma F}_{5^2}$
is the Galois closure of $K_0/k$.
We have that $T$ and $T^2+T+1$ are irreducible in ${\ma F}_5(T)$
since $\zeta_3\notin {\ma F}_5$ and $T^2+T+1=\frac{T^3-1}{T-1}=
(T-\zeta_3)(T-\zeta_3^2)$. In fact ${\ma F}_5(\zeta_3)={\ma F}_{5^2}
={\ma F}_{25}$. 

Let ${\mathcal P}_T$ and ${\mathcal P}_{
T^2+T+1}$ be the prime divisors in $k={\ma F}_5(T)$ corresponding
to $T$ and $T^2+T+1$ respectively. 
The infinite prime $\p$ is unramified in $K_0/k$ and in $K/k$ because $\deg
D(T)=3$. Let $t_0(K_0)$ and $t_0(K)$ be given by (\ref{Eq3.2''}) with
respect to the fields $K_0$ and $K$ respectively. Since $\gamma=1$, from
(\ref{Eq5.1*}) we obtain
\begin{gather*}
t_0(K_0)=\gcd_{0\leq i\leq 2}\big\{\big[{\ma F}_5(\zeta_3^i):{\ma F}_5\big]\big\}=
\gcd\{1,2,2\}=1
\intertext{and}
t_0(K)=\gcd_{0\leq i\leq 2}\big\{\big[{\ma F}_{5^2}(\zeta_3^i):{\ma F}_5\big]\big\}
=\gcd_{0\leq i\leq 2}\{[{\ma F}_{5^2}:{\ma F}_5]\}= \gcd\{2,2,2\}=2.
\end{gather*}
Since $t_0(K_0)\neq t_0(K)$, it follows that $\g{(K_0)}\neq \g K$.

Let $F_0=F_{{\mathcal P}_T} F_{{\mathcal P}_{T^2+T+1}}$. We have
\begin{gather*}
[F_{{\mathcal P}_T}:k]=c_{{\mathcal P}_T}=\gcd(e_{{\mathcal P}_T},
5^{d_{{\mathcal P}_T}}-1)=(3,4)=1.
\intertext{Therefore $F_{{\mathcal P}_T}=k$.
Now}
[F_{{\mathcal P}_{T^2+T+1}}:k]=c_{{\mathcal P}_{T^2+T+1}}=
\gcd(e_{{\mathcal P}_{T^2+T+1}}, q^{\deg {\mathcal P}_{T^2+T+1}}-1)=
\gcd(3,24)=3.
\end{gather*}

Since $q-1=e_{\p}(k(\Lambda_{T^2+T+1})|k)=4$ it
follows that $F_{{\mathcal P}_{T^2+T+1}}\subseteq k(\Lambda_{T^2
+T+1})^+$ and $F_0=F_{{\mathcal P}_{T^2+T+1}}=F_0\cap
k(\Lambda_{T(T^2+T+1)})^+$. Hence $F_0$ is the
unique subfield of $k(\Lambda_{T^2+T+1})$ of degree $3$ 
over $k$ and $F_0=F$.

We have that $F_{{\mathcal P}_{T^2+T+1}}\cdot {\ma F}_{5^2}/
{\ma F}_{5^2}(T)$ is a Kummer extension of degree $3$ where
the finite ramified primes are the primes in ${\ma F}_{5^2}(T)$
dividing $T^2+T+1=(T-\zeta_3)(T-\zeta_3^2)$. Since $\p$ decomposes
fully in $F_{{\mathcal P}_{T^2+T+1}}/k$, ${\eu p}_{\infty}$, the infinite prime
in ${\ma F}_{5^2}(T)$, decomposes fully in $F_{{\mathcal P}_{T^2+T+1}}
\cdot {\ma F}_{5^2}/{\ma F}_{5^2}(T)$. Therefore 
$F_{{\mathcal P}_{T^2+T+1}} \cdot {\ma F}_{5^2}={\ma F}_{5^2}(T)
\big(\sqrt[3]{(-1)^{\deg Q}Q(T)}\big) = {\ma F}_{5^2}(T)
\big(\sqrt[3]{Q(T)}\big)$ with $\deg Q(T)=3$ and $T-\zeta_3, T-\zeta_3^2$ are the
unique irreducible polynomials dividing $Q(T)$.
\[
\xymatrix{
F_{{\mathcal P}_{T^2+T+1}}\ar@{-}[r]
\ar@{-}[d]_3 & F_{{\mathcal P}_{T^2+T+1}}
\cdot {\ma F}_{5^2}={\ma F}_{25}(T)\big(\sqrt[3]{Q(T)}\big)
\ar@{-}[d]^3\\
k\ar@{-}[r]&{\ma F}_{25}(T)
}
\]

It follows that 
\begin{gather*}
F_{{\mathcal P}_{T^2+T+1}}
\cdot {\ma F}_{5^2}={\ma F}_{25}(T)\Big(\sqrt[3]{(T-\zeta_3)(T-\zeta_3^2)^2}\Big)
={\ma F}_{25}(T)\Big(\sqrt[3]{(T-\zeta_3)^2(T-\zeta_3^2)}\Big).
\intertext{Now, since $F_0=F$, from Remark \ref{R3.5} we obtain}
\g K=KF{\ma F}_{5^2}={\ma F}_{25}\Big(T,\sqrt[3]{T(T-\zeta_3)(T-\zeta_3^2)},
\sqrt[3]{(T-\zeta_3)(T-\zeta_3^2)^2}\Big).
\end{gather*}

Let $\g {K^{\prime}}$ be the genus
field of $K/{\ma F}_{25}(T)$. We may apply Peng's
Theorem. With the notations from \cite[Theorem 5.2]{MaRzVi2015}, we have
$r=3, P_1=T, P_2=T-\zeta_3, P_3=
T-\zeta_3^2, \gamma =1, \alpha=(-1)^{\deg D}\gamma =-1\in({\ma F}_{25}^{
\ast})^3, a_1=a_2=2$. Thus
\[
\g {K^{\prime}}={\ma F}_{25}\Big(T, \sqrt[3]{T(T-\zeta_3^2)^2},\sqrt[3]
{(T-\zeta_3)(T-\zeta_3^2)^2}\Big).
\]

We also have 
\begin{multline*}
{\ma F}_{25}\Big(T,\sqrt[3]{T(T-\zeta_3)(T-\zeta_3^2)},
\sqrt[3]{(T-\zeta_3)(T-\zeta_3^2)^2}\Big) \\
={\ma F}_{25}\Big(T, \sqrt[3]{T(T-\zeta_3^2)^2},\sqrt[3]
{(T-\zeta_3)(T-\zeta_3^2)^2}\Big).
\end{multline*}
Therefore $\g K=\g {K^{\prime}}$ (see Remark \ref{R5.4}).

Finally, from Remark \ref{R3.5} $\g {(K_0)}=KF=k\big(T,\sqrt[3]{D(T)}\big)
F_{{\mathcal P}_{T^2+T+1}}$.
}}
\end{example}

\begin{remark}\label{R5.4}{\rm{
Let $k={\ma F}_q(T)$ and let $k_n$ be the extension of constants 
of $k$ of degree $n$.
Let $K$ be any finite extension of $k$ such that ${\ma F}_{q^n}
\subseteq K$. Let $\g K$ and
$\g {K^{\prime}}$ be the genus fields of $K/k$ and $K/k_n$ respectively.
Since the infinite prime divisors of $K$ decompose fully in $\g K$
and in $\g {K^{\prime}}$ and $\g K/K$ and $\g {K^{\prime}}/K$ are abelian
and unramified, it follows that $\g K=\g {K^{\prime}}$.
}}
\end{remark}

\end{document}